\documentstyle[amscd,amssymb,verbatim,epsf]{amsart}
\begin{document}

\theoremstyle{plain}
\newtheorem{Thm}{Theorem}
\newtheorem{Cor}{Corollary}
\newtheorem{Con}{Conjecture}
\newtheorem{Main}{Main Theorem}

\newtheorem{Lem}{Lemma}
\newtheorem{Prop}{Proposition}

\theoremstyle{definition}
\newtheorem{Def}{Definition}
\newtheorem{Note}{Note}

\theoremstyle{remark}
\newtheorem{notation}{Notation}
%\textbf{R}newcommand{\thenotation}{}

\errorcontextlines=0
\numberwithin{equation}{section}
%\textbf{R}newcommand{\rm}{\normalshape}%

\title[$\lambda$-statistically quasi-Cauchy sequences]%
   {$\lambda$-statistically quasi-Cauchy sequences}
\author{H\"usey\.{I}n \c{C}akalli*, Ayse Sonmez** , and Cigdem Gunduz Aras*** \\
*Maltepe University, Marmara E\u{g}\.{I}t\.{I}m K\"oy\"u, TR 34857, \.{I}stanbul-Turkey \; \; \; \; \; Phone:(+90216)6261050 ext:2248, \;  fax:(+90216)6261113 \\ **Gebze Institute of Technology, Department of Mathematics,  Gebze-Kocaeli, Turkey Phone: (+90262)6051389 \\*** Kocaeli University, Department of Mathematics, Kocaeli, Turkey Phone:(+90262)3032102}

\address{H\"usey\.{I}n \c{C}akall\i \\
          Maltepe University, Department of Mathematics, Marmara E\u{g}\.{I}t\.{I}m K\"oy\"u, TR 34857, Maltepe, \.{I}stanbul-Turkey \; \; \; \; \; Phone:(+90216)6261050 ext:2248, \;  fax:(+90216)6261113}

\email{hcakalli@@maltepe.edu.tr; hcakalli@@gmail.com}

\address{Ay\c{s}e S\"onmez\\
                     Department of Mathematics, Gebze Institute of Technology, Cayirova Campus 41400 Gebze- Kocaeli, Turkey Phone: (+90262)6051389}

\email{asonmez@@gyte.edu.tr; ayse.sonmz@@gmail.com}

\address{Cigdem Gunduz Aras\\
                    Kocaeli University, Department of Mathematics, Kocaeli, Turkey Phone:(+90262)3032102}

\email{caras@@kocaeli.edu.tr; carasgunduz@@gmail.com }

\date{\today}

\maketitle

\begin{abstract}
The main object of this paper is to investigate $\lambda$-statistically quasi-Cauchy sequences. A real valued function $f$ defined on a subset $E$ of $\textbf{R}$, the set of real numbers, is called $\lambda$-statistically ward continuous on $E$ if it preserves $\lambda$-statistically quasi-Cauchy sequences of points in $E$. It turns out that uniform continuity coincides with $\lambda$-statistically ward continuity on $\lambda$-statistically ward compact subsets. 
\end{abstract}

\maketitle

\section{Introduction}

A function $f:\textbf{R} \longrightarrow \textbf{R}$ is continuous if and only if it preserves Cauchy sequences.  Using the idea of continuity of a real function in terms of sequences, many kinds of continuities were introduced and investigated, not all but some of them we recall in the following: slowly oscillating continuity  (\cite{CakalliSlowlyoscillatingcontinuity}), quasi-slowly oscillating continuity (\cite{DikandCanak}), ward continuity (\cite{CakalliForwardcontinuity}), $\delta$-ward continuity (\cite{CakalliDeltaquasiCauchysequences}), statistical ward continuity,  (\cite{CakalliStatisticalwardcontinuity}), and $N_{\theta}$-ward continuity (\cite{CakalliNthetawardcontinuity}) which enabled some authors to obtain some characterizations of uniform continuity in terms of sequences in the sense that a function preserves either quasi-Cauchy sequences or slowly oscillating sequences (\cite{Vallin}, and \cite{BurtonandColemanQuasiCauchySequences}). In \cite{MursaleenLamdastatisticalconvergence}, the concept of statistical convergence was generalized to $\lambda$-statistically convergence.

In this paper, we investigate the concept of $\lambda$-statistically ward continuity, and examine its properties.

\maketitle

\section{Preliminaries}
\normalfont{}
Boldface letters $\boldsymbol{\alpha}$, $\bf{x}$, $\bf{y}$, $\bf{z}$, ... will be used for sequences $\boldsymbol{\alpha}=(\alpha_{n})$, $\textbf{x}=(x_{n})$, $\textbf{y}=(y_{n})$, $\textbf{z}=(z_{n})$, ... of points in $\textbf{R}$ for the sake of abbreviation. $s$ and $c$ will denote the set of all sequences, and the set of convergent sequences of points in \textbf{R}.

A sequence $(\alpha _{n})$ of points in $\textbf{R}$ is quasi-Cauchy if $(\Delta \alpha _{n})$ is a null sequence, where $\Delta \alpha _{n}=\alpha _{n+1}-\alpha _{n}$. These sequences were named as quasi-Cauchy by Burton and Coleman \cite[page 328]{BurtonandColemanQuasiCauchySequences}, while they  were called as forward convergent to $0$ sequences in \cite[page 226]{CakalliForwardcontinuity}. Quasi-Cauchy sequences arise in diverse situations, and it is often difficult to determine whether or not they converge, and if so, to which limit. It is easy to construct a zero-one sequence such that the quasi-Cauchy average sequence does not converge. The usual constructions have a somewhat artificial feeling. Nevertheless, there are sequences which seem natural, have the quasi-Cauchy property, and do not converge. On the other hand, the sequence of averages of $0$ s and $1$ s is always a quasi-Cauchy sequence: let $\textbf{x}:= (x_{n})$ be a sequence such that for each nonnegative integer $n$, $x_{n}$ is either $0$ or $1$. For each positive integer $n$ set  $a_{n}=\frac{x_{1}+x_{2}+...+x_{n}}{n}$ . Then $a_{n}$ is the arithmetic mean average of the sequence up to time or position $n$. Clearly
for each $n$, $0 \leq a_{n} \leq 1$, and  $|a_{n+1}-a_{n}|  \leq {\frac{1}{n}}$. Thus $(a_{n})$ is a quasi-Cauchy sequence.

The concept of statistical convergence is a generalization of the usual notion of convergence that, for real-valued sequences, parallels the usual theory of convergence (see \cite{Fridy}).

Let  $\boldsymbol{\lambda} = (\lambda_{n})$ be a non-decreasing sequence of positive numbers tending to $\infty$ such that $\lambda_{n+1}\leq \lambda_{n} + 1$, $\lambda_{1} = 1$. The generalized de la Valee-Pousin mean of a sequence $\boldsymbol{\alpha} = (\alpha_{k})$ is defined by
$$t_{n}(\boldsymbol{\alpha}):=\frac{1}{\lambda_{n}} \sum^{}_{k\in{I_{n}}} \alpha_{k}$$
where $I_{n} = [n - \lambda_{n}+1, n]$. A sequence $\boldsymbol{\alpha} = (\alpha_{k})$ is said to be $(V, \lambda)$-summable to a number $L$ if
\[t_{n}(\boldsymbol{\alpha})\longrightarrow L\; \; \;   as\;  n\longrightarrow\infty,\]
%This is denoted by $V_{\lambda}-lim x_{k}=L$.
a sequence $\boldsymbol{\alpha} = (\alpha_{k})$ is said to be $[V, \lambda]$-summable to a number $L$ or strongly $(V, \lambda)$-summable to a number $L$  (see \cite{LeindlerVallePousinche}) if
$$\lim_{n\rightarrow\infty}\frac{1}{\lambda_{n} }\sum^{}_{k\in{I_{n}}} |\alpha_{k}-L|=0,$$
and a sequence $\boldsymbol{\alpha} = (\alpha_{k})$ of points in $\textbf{R}$ is called to be $\lambda$-statistically convergent to an element $L$ of $\textbf{R}$ if
\[
\lim_{n\rightarrow\infty}\frac{1}{\lambda_{n} }|\{k\in I_{n}: |\alpha_{k}-L| \geq\varepsilon\}|=0
\]
for every positive real number $\varepsilon$ (\cite{MursaleenLamdastatisticalconvergence}). This is denoted by $S_{\lambda}-lim \alpha_{k}=\alpha_{0}$.
If we write $\lambda_{n}=n$ for all $n\in{\textbf{N}}$, then we get statistical convergence. Throughout this paper, $S$, $S_{\lambda}$ will denote the set of statistically convergent and $\lambda$-statistically convergent sequences of points in $\textbf{R}$, respectively.

\maketitle
\section{$\lambda$-statistically ward continuity}
Connor and Grosse-Erdman (\cite{ConnorandGrosse}) gave sequential definitions of continuity for real functions calling $G$-continuity instead of $A$-continuity and their results covers the earlier works related to $A$-continuity where a method of sequential convergence, or briefly a method, is a linear function $G$ defined on a linear subspace of $s$, denoted by $c_{G}$, into $\textbf{R}$. A sequence $\textbf{x}=(x_{n})$ is said to be $G$-convergent to $L$ if $\textbf{x}\in c_{G}$ and $G(\textbf{x})=L$. A method $G$ is called regular if every convergent sequence $\textbf{x}=(x_{n})$ is $G$-convergent with $G(\textbf{x})=\lim \textbf{x}$. A method $G$ is called subsequential if whenever $\textbf{x}$ is $G$-convergent with $G(\textbf{x})=L$, then there is a subsequence $(x_{n_{k}})$ of $\textbf{x}$ with $\lim_{k} x_{n_{k}}=L$. In particular, $\lim$ denotes the limit function $\lim \textbf{x}=\lim_{n}x_{n}$ on the linear space $c$, $st-\lim$ denotes the statistical limit function $st-\lim \textbf{x}=st-\lim_{n}x_{n}$ on the linear space $S$, $S_{\theta}-\lim$ denotes the lacunary statistical limit function $S_{\theta}-\lim \textbf{x}=S_{\theta}-\lim_{n}x_{n}$ on the linear space $S_{\theta}$,  $V_{\lambda}-\lim$ denotes the  $V_{\lambda}$-statistical limit function $V_{\lambda}-\lim \textbf{x}=V_{\lambda}-\lim_{n}x_{n}$ on the linear space $V_{\lambda}$,  $[V, \lambda]-\lim$ denotes the  $[V, \lambda]$-limit function (or strongly $(V,\lambda)$-limit function) $[V, \lambda]-\lim \textbf{x}=[V, \lambda]-\lim_{n}x_{n}$ on the linear space $[V, \lambda]$, and $S_{\lambda}-\lim$ denotes the  $S_{\lambda}$-statistical limit function $S_{\lambda}-\lim \textbf{x}=S_{\lambda}-\lim_{n}x_{n}$ on the linear space $S_{\lambda}$. All sequential methods exempt lacunary statistical method mentioned a few lines above are all regular without any restriction, however lacunary statistical method is regular under the assumption that $\lim inf_{r}\; q_{r}>1$ where a real sequence $(x_{k})$ is called lacunary statistically convergent to an element $L$ of $\textbf{R}$ if
\[
\lim_{r\rightarrow\infty}\frac{1}{h_{r}} |\{k\in I_{r}: |x_{k}-L|\geq{\epsilon} \}|=0,
\]
for every $\epsilon>0$, $I_{r}=(k_{r-1},k_{r}]$ and $k_{0}=0$, $h_{r}:k_{r}-k_{r-1}\rightarrow
\infty$ as $r\rightarrow\infty$ and $\theta=(k_{r})$ is an increasing sequence of positive integers. In the sequel, we will always assume that $\lim inf_{r}\; q_{r}>1$. $S_{\lambda}$-statistical sequential method is subsequential as Statistical sequential method and lacunary statistical sequential method are subsequential, so is $S_{\lambda}$-statistical sequential method.

A subset $E$ of $\textbf{R}$ is called $G$-sequentially compact if any sequence $\textbf{x}$ of points in $E$ has a $G$-sequentially convergent subsequence $\textbf{z}=(x_{n_{k}})$ such that $G(\textbf{z})\in {E}$. A subset $E$ of $\textbf{R}$ is $\lambda$-statistically sequentially compact if any sequence $\textbf{x}$ of points in $E$ has a $\lambda$-statistically convergent subsequence whose $\lambda$-statistical limit is in $E$. We see that this is a special case of $G$-sequential compactness where $G=S_{\lambda}-lim$.
$\lambda$-statistically sequentially compactness of a subset $E$ of $\textbf{R}$ coincides not only with ordinary (sequential) compactness, but also statistically sequentially compactness, and lacunary statistically sequentially compactness.
 A function $f$ is called $G$-continuous at a point $u$ provided that whenever a sequence $\textbf{x}=(x_{n})$ of terms in the domain of $f$ is $G$-convergent to $u$, then the sequence $f(\textbf{x})=(f(x_{n}))$ is $G$-convergent to $f(u)$. Writing $G=S_{\lambda}$, we get $S_{\lambda}$-sequential continuity or $\lambda$-statistically sequential continuity, explicitly we say that a real valued function $f$ defined on a subset $E$ of $R$ is called $\lambda$-statistically sequentially continuous at a point $\alpha_{0}$ if it preserves $\lambda$-statistically convergent sequences, i.e. $(f(\alpha_{k}))$ is a $\lambda$-statistically convergent to $f(\alpha_{0})$ whenever $(\alpha_{k})$ is $\lambda$-statistically convergent to $\alpha_{0}$. $\lambda$-statistically sequentially continuity of a real valued function defined on a subset of $\textbf{R}$ coincides with not only ordinary (sequential) continuity, but also each one of the continuities, statistically sequentially continuity, and lacunary statistically sequentially continuity.

A sequence $(\alpha_{k})$ of points in $\textbf{R}$ is called to be $\lambda$-statistically quasi-Cauchy if $S_{\lambda}-lim \Delta\alpha_{k}=0$.
Now we give the following interesting examples which show emphasis the interest in different research areas.

\textbf{Example 1.} (\cite{WinklerMathematicalPuzzles}) Let $n$ be a positive integer. In a group of $n$ people, each person selects at random and simultaneously another person of the group. All of the selected persons are then removed from the group, leaving a random number $n_{1} < n$ of people which form a new group. The new group then repeats independently the selection and removal thus described, leaving $n_{2} < n_{1}$ persons, and so forth until either one person remains, or no persons remain. Denote by $p_n$ the probability that, at the end of this iteration initiated with a group of $n$ persons, one person remains. Then the sequence
$\textbf{p} = (p_{1}, p_{2}, · · ·, p_{n},...)$  is a $\lambda$-statistically quasi-Cauchy sequence, and $lim p_n$ does not exist.

\textbf{Example 2.} (\cite{KeaneUnderstandingErgodicity}) Let $n$ be a positive integer. In a group of $n$ people, each person selects independently and at random one of three subgroups to which to belong, resulting in three groups with random numbers $n_{1}$, $n_{2}$, $n_{3}$ of members; $n_{1} + n_{2} + n_{3} = n$. Each of the subgroups is then partitioned independently in the same manner to form three sub subgroups, and so forth. Subgroups having no members or having only one member are removed from the process. Denote by $t_{n}$ the expected value of the number of iterations up to complete removal, starting initially with a group of $n$ people. Then the sequence $(t_{1}, \frac{t_{2}}{2}, \frac{t_{3}}{3},...,\frac{t_{k}}{k},...)$ is a bounded nonconvergent $\lambda$-statistically quasi-Cauchy sequence.

Now we state the definition of $\lambda$-statistically ward continuity in the following.
\begin{Def}
A real valued function $f$ defined on a subset $E$ of $\textbf{R}$ is called $\lambda$-statistically ward continuous on $E$ if it preserves $\lambda$-statistically quasi-Cauchy sequences of points in $E$, i.e. $(f(\alpha_{k}))$ is a $\lambda$-statistically quasi-Cauchy sequence whenever $(\alpha_{k})$ is a $\lambda$-statistically quasi-Cauchy sequences of points in $E$.
\end{Def}

We note that this definition of continuity cannot be obtained by any $A$-continuity, i.e., by any summability matrix $A$, even by the summability matrix $A=(a_{nk})$ defined by
$$a_{nk}=\frac{1}{\lambda_{n}}\; for\; k=n+1, \;\;\;and\;\;\; a_{nk}=-\frac{1}{\lambda_{n}}\; for\; k=n.$$
However, for this special summability matrix $A$, if $A$-continuity of $f$ at the point $0$ implies $\lambda$-statistically ward continuity of $f$, then $f(0)=0$; and if $\lambda$-statistically ward continuity of $f$ implies $A$-continuity of $f$ at the point $0$, then $f(0)=0$.

Sum of two $\lambda$-statistically ward continuous functions is $\lambda$-statistically ward continuous, but product of $\lambda$-statistically ward continuous functions need not be $\lambda$-statistically ward continuous.

We give the following theorem the proof of which also can be obtained by considering \cite[Theorem 9]{CakalliandHazarikaIdealquasiCauchysequences}.

\begin{Thm} \label{lacunarystatisticalwardcontinuityimplieslacunarystatiscontinuity} If a real valued function is $\lambda$-statistically ward continuous on a subset $E$ of $\textbf{R}$, then it is $\lambda$-statistically sequentially continuous on $E$.
\end{Thm}
\begin{pf} Suppose that $f$ is a $\lambda$-statistically ward continuous function on a subset $E$ of $\textbf{R}$. Let $(x_{n})$ be a $\lambda$-statistically quasi-Cauchy sequence of points in $E$. Then the sequence $$(x_1, x_0, x_2, x_0, x_3, x_0,..., x_{n-1}, x_0, x_n, x_0,...)$$ is a $\lambda$-statistically quasi-Cauchy sequence. Since $f$ is $\lambda$-statistically ward continuous, the sequence $$(y_n)=(f(x_{1}),f(x_{0}),f(x_{2}),f(x_{0}),...,f(x_{n}),f(x_{0}),...)$$ is a  $\lambda$-statistically quasi-Cauchy sequence. Therefore $S_{\lambda}-\lim_{n\rightarrow \infty} \Delta y_n=0$. Hence  $S_{\lambda}-\lim_{n\rightarrow \infty} [f(x_{n})-f(x_{0})]=0$. It follows that the sequence $(f(x_{n}))$ is $\lambda$-statistically summable to $f(x_{0})$. This completes the proof of the theorem.
\end{pf}

\begin{Cor} \label{lacunarystatisticalwardcontinuityimpliesordinarycontinuity} If a real valued function is $\lambda$-statistically ward continuous on a subset $E$ of \textbf{R}, then it is continuous on $E$.
\end{Cor}

\begin{pf}
The proof follows immediately from the preceding theorem so is omitted.
\end{pf}

Now we prove the following theorem.
\begin{Thm} \label{TheoremuniformlycontinuousfunctiononEsendsquasiCauchytolacunarystatisticallquasiCauchy}
If a real valued function $f$ is uniformly continuous on a subset $E$ of $\textbf{R}$, then $(f(x_{n}))$ is $\lambda$-statistically quasi-Cauchy whenever $(x_{n})$ is a quasi-Cauchy sequence of points in $E$.
\end{Thm}
\begin{pf} Let $f$ be uniformly continuous on $E$. Take any quasi-Cauchy sequence $(x_{n})$ of points in $E$. Let $\varepsilon$ be any positive real number. Since $f$ is uniformly continuous, there exists a $\delta>0$ such that $|f(x)-f(y)|<\varepsilon$ whenever $|x-y|<\delta$. As $(x_{n})$ is a quasi-Cauchy sequence, for this $\delta$ there exists an $n_{0}\in{\textbf{N}}$ such that $|x_{n+1}-x_{n}|<\delta$ for $n\geq n_{0}$. Therefore $|f(x_{n+1})-f(x_{n})|<\varepsilon$ for $n\geq n_{0}$, so the number of indices $k$ for which $|f(x_{n+1})-f(x_{n})|\geq\varepsilon$ is less than $n_{0}$. Hence

$\lim_{r\rightarrow\infty}\frac{1}{\lambda_{r} }|\{k\in I_{r}:|f(x_{n+1})-f(x_{n})|\geq\varepsilon\}|\leq \lim_{r\rightarrow\infty}\frac{n_{0}}{\lambda_{r} }=0$.\\
This completes the proof of the theorem.
\end{pf}

On the other hand, any continuous function on a compact subset $E$ of $\textbf{R}$ is uniformly continuous on $E$. It is also true for a regular subsequential method $G$ that any $\lambda$-statistically ward continuous function on a $G$-sequentially compact subset $E$ of $\textbf{R}$ is also uniformly continuous on $E$ (see \cite{CakalliSequentialdefinitionsofcompactness}). Furthermore, for $\lambda$-statistically ward continuous functions defined on a $\lambda$-statistically ward compact subset of \textbf{R}, we have the following.

\begin{Thm} \label{Lamdastatisticalwardcontinuousfunctiononlacunarystatiswardcompactsubstisunifotmlycontinuous} Any $\lambda$-statistically ward continuous real valued function on a $\lambda$-statistically ward compact subset of $\textbf{R}$ is uniformly continuous.

\end{Thm}
\begin{pf} Let $E$ be a $\lambda$-statistically ward compact subset $E$ of $\textbf{R}$ and let $f:E\longrightarrow$ $\textbf{R}$ be a  $\lambda$-statistically ward continuous function on $E$. Suppose that $f$ is not uniformly continuous on $E$ so that there exists an  $\varepsilon_{0} > 0$ such that for any $\delta >0$, there are $x, y \in{E}$ with $|x-y|<\delta$ but $|f(x)-f(y)| \geq \varepsilon_{0}$. For each positive integer $n$, there exist $\alpha_{n}$ and $\beta_{n}$ such that $|\alpha _{n}-\beta_{n}|<\frac{1}{\lambda_{n}}$, and $|f(\alpha _{n})-f(\beta_{n})|\geq \varepsilon_{0}$. Since $E$ is $\lambda$-statistically ward compact, there exists a $\lambda$-statistically quasi-Cauchy subsequence $(\alpha _{n_{k}})$ of the sequence $(\alpha _{n})$. It is clear that the corresponding subsequence $(\beta_{n_{k}})$ of the sequence $(\beta_{n})$ is also $\lambda$-statistically quasi-Cauchy, since $(\beta_{n_{k+1}}-\beta_{n_{k}})$ is $\lambda$-statistically convergent to $0$ which follows from the following lines: for each $\varepsilon$,
 $$\{k\in I_{r}: |\beta_{n_{k+1}}-\beta_{n_{k}}| \geq \varepsilon\} \subset {\{k\in I_{r}: |\beta_{n_{k+1}}-\alpha_{n_{k+1}}|\geq \frac{\varepsilon}{3}\} \cup \{k\in I_{r}: |\alpha_{n_{k+1}}-\alpha_{n_{k}}| \geq  \frac{\varepsilon}{3} \}}$$ $$ \cup \{k\in I_{r}: |\alpha_{n_{k}}-\beta_{n_{k}}| \geq \frac{\varepsilon}{3}\}.$$ It follows from this inclusion that \\
$|\{k\in I_{r}: |\beta_{n_{k+1}}-\beta_{n_{k}}| \geq \varepsilon\}|$
$\leq {|\{k\in I_{r}: |\beta_{n_{k+1}}-\alpha_{n_{k+1}}|\geq \frac{\varepsilon}{3}\}}|$\\ $+
|\{k\in I_{r}: |\alpha_{n_{k+1}}-\alpha_{n_{k}}| \geq \frac{\varepsilon}{3}\}|+$
$|\{k\in I_{r}: |\alpha_{n_{k}}-\beta_{n_{k}}| \geq \frac{\varepsilon}{3}\}|$.

Hence

$\lim_{r\rightarrow\infty}\frac{1}{\lambda_{r} }|\{k\in I_{r}: |\beta_{n_{k+1}}-\beta_{n_{k}}| \geq \varepsilon\}|$
$\leq {\lim_{r\rightarrow\infty}\frac{1}{\lambda_{r} }|\{k\in I_{r}: |\beta_{n_{k+1}}-\alpha_{n_{k+1}}|\geq \frac{\varepsilon}{3} \}}|+$
$\lim_{r\rightarrow\infty}\frac{1}{\lambda_{r} }|\{k\in I_{r}: |\alpha_{n_{k+1}}-\alpha_{n_{k}}| \geq \frac{\varepsilon}{3}\}|+$$\lim_{r\rightarrow\infty}\frac{1}{\lambda_{r} }|\{k\in I_{r}:|\alpha_{n_{k}}-\beta_{n_{k}}| \geq \frac{\varepsilon}{3}\}|=0$

On the other hand, it follows from the equality $$\alpha_{n_{k+1}}-\beta_{n_{k}}=\alpha_{n_{k+1}}-\alpha_{n_{k}}+\alpha_{n_{k}}-\beta_{n_{k}}$$ that the sequence $(\alpha_{n_{k+1}}-\beta_{n_{k}})$ is $\lambda$-statistically convergent to $0$.
Hence the sequence
$$(a_{n_{1}}, \beta_{n_{1}}, \alpha _{n_{2}}, \beta_{n_{2}}, \alpha _{n_{3}}, \beta_{n_{3}},..., \alpha _{n_{k}}, \beta_{n_{k}},...)$$
is $\lambda$-statistically quasi-Cauchy.
But the transformed sequence
$$(f(\alpha _{n_{1}}), f(\beta_{n_{1}}), f(\alpha _{n_{2}}), f(\beta_{n_{2}}), f(\alpha _{n_{3}}), f(\beta_{n_{3}}),..., f(\alpha _{n_{k}}), f(\beta_{n_{k}}),...)$$
is not $\lambda$-statistically quasi-Cauchy. Thus $f$ does not preserve $\lambda$-statistically quasi-Cauchy sequences. This contradiction completes the proof of the theorem.
\end{pf}

\begin{Cor} \label{Lamdatatiswardcontinuousfunctiononboundedsubsetisunifotmlycontinuous} If a real valued function is $\lambda$-statistically ward continuous on a bounded subset $E$ of $\textbf{R}$, then it is uniformly continuous on $E$.
\end{Cor}

\begin{pf}
The proof follows from the fact that any bounded subset of $\textbf{R}$ is $\lambda$-statistically ward compact.
%Theorem 3 in \cite[page 1622]{CakalliStatisticalquasiCauchysequences}.
\end{pf}

\begin{Cor} \label{LamdatatiswardcontinuousfunctiononNthetawardcompactsubsetisunifotmlycontinuous} If a real valued function is $\lambda$-statistically ward continuous on an $N_{\theta}$-ward compact subset $E$ of \textbf{R}, then it is uniformly continuous on $E$.
\end{Cor}

\begin{pf}
The proof follows from \ref{Lamdastatisticalwardcontinuousfunctiononlacunarystatiswardcompactsubstisunifotmlycontinuous} and \cite[Theorem 3.3]{CakalliNthetawardcontinuity}.
\end{pf}

We give below that any real $\lambda$ statistically-ward continuous function defined on an interval is uniformly continuous. First we give the following lemma.

\begin{Lem} \label{Lemmaanypairofsequencehasalamdastatisticallystatisticallyquasi}
If $(\xi_{n}, \eta_{n})$ is a sequence of ordered pairs of points in an interval such that $\lim_{n\rightarrow\infty} |\xi_{n}-\eta_{n}|=0$, then there exists a  $\lambda$-statistically quasi-Cauchy sequence $(\alpha_{n})$ with the property that for any positive integer $i$ there exists a positive integer $j$ such that $(\xi_{i}, \eta_{i})=(\alpha_{j-1}, \alpha_{j})$.
\end{Lem}

\begin{Thm} \label{TheoremLamdastatisticallywardcontinuityonanintervalimpliesuniformlycontinuity}
If a real valued function defined on an interval $E$ is $\lambda$-statistically ward continuous, then it is uniformly continuous.
\end{Thm}
\begin{pf}
Suppose that $f$ is not uniformly continuous on $E$. Then there is an $\varepsilon_{0} > 0$ such that for any $\delta > 0$ there exist $x, y \in {E}$ with $|x - y| < \delta$ but $| f (x) - f (y)|\geq {\varepsilon_{0}}$. For every $n \in{\textbf{N}}$ fix $\xi_{n}$, $\eta_{n}\in{E}$ with $|\xi_{n} - \eta_{n}| < \frac{1}{n}$ and $| f (\xi_{n}) - f (\eta_{n})| \geq \varepsilon_{0}$. By Lemma \ref{Lemmaanypairofsequencehasalamdastatisticallystatisticallyquasi}, there exists a $\lambda$- statistically quasi-Cauchy sequence $(\alpha_{i})$ such that for any integer $i \geq {1}$ there exists a $j$ with $\xi_{i} = \alpha_{j}$ and $\eta_{i} = \alpha_{j+1}$. This implies that $| f (\alpha_{j+1} ) - f (\alpha_{j})| \geq \varepsilon_{0}$; hence $( f (\alpha_{i} ))$ is not  $\lambda$-statistically quasi-Cauchy. Thus $f$ does not preserve $\lambda$-statistically quasi-Cauchy sequences. This completes the proof of the theorem.
\end{pf}

Since the sequence constructed in Lemma \ref{Lemmaanypairofsequencehasalamdastatisticallystatisticallyquasi} is also quasi-Cauchy, we see that the statement $(f(x_{n}))$ is $\lambda$ statistically quasi-Cauchy whenever $(x_{n})$ is quasi-Cauchy sequence of points in $E$ implies the uniform continuity of $f$ on $E$. Now combining Theorem \ref{TheoremuniformlycontinuousfunctiononEsendsquasiCauchytolacunarystatisticallquasiCauchy} with this observation we have the following result.

\begin{Cor}
Let $f$ be a real valued function defined on an interval $E$. Then  $f$ is uniformly continuous on $E$ if and only if $(f(x_{n}))$ is $\lambda$-statistically quasi-Cauchy whenever $(x_{n})$ is quasi-Cauchy sequence of points in $E$.
\end{Cor}

\begin{Cor}
Let $f$ be a real valued function defined on an interval $E$. Then the following statements are equivalent:\\
(a)   if $(f(x_{n}))$ is $\lambda$-statistically quasi-Cauchy whenever $(x_{n})$ is quasi-Cauchy sequence of points in $E$.\\
(b)   if $(f(x_{n}))$ is $N_{\theta}$ quasi-Cauchy whenever $(x_{n})$ is quasi-Cauchy sequence of points in $E$.
\end{Cor}

\begin{pf}
The proof follows from Theorem \ref{TheoremLamdastatisticallywardcontinuityonanintervalimpliesuniformlycontinuity}, and \cite[Theorem 1 and Theorem 2]{CakalliandKaplanAstudyonNthetawardcontinuity} so is omitted.
\end{pf}

\begin{Cor}
If a real valued function defined on an interval is $\lambda$-statistically ward continuous, then it is ward continuous.
\end{Cor}

\begin{pf}
The proof follows from Theorem \ref{TheoremLamdastatisticallywardcontinuityonanintervalimpliesuniformlycontinuity}, so it is omitted.
\end{pf}
\begin{Cor}
If a real valued function defined on an interval is $\lambda$-statistically ward continuous, then it is slowly oscillating continuous.
\end{Cor}
\begin{pf}
The proof follows from Theorem \ref{TheoremLamdastatisticallywardcontinuityonanintervalimpliesuniformlycontinuity}, so it is omitted.
\end{pf}

\section{Conclusion}
In this paper, the concept of $\lambda$-statistically ward continuity of a real function is
investigated. In this investigation we have obtained theorems related to $\lambda$-statistically ward continuity, $N_{\theta}$-ward continuity, statistically ward continuity, forward continuity, uniform continuity, $\lambda$-statistically ward compactness, $N_{\theta}$-ward compactness, statistically ward compactness, and compactness.
For a further study, we suggest to investigate $\lambda$-statistically quasi-Cauchy sequences of fuzzy points, and $\lambda$-statistically ward continuity for the fuzzy functions (see \cite{SavasOnasymptoticallylacunarystatisticalequivalentsequencesoffuzzynumbers} for the definitions and  related concepts in fuzzy setting). However due to the change in settings, the definitions and methods of proofs will not always be analogous to those of the present work. For another further study we suggest to investigate $\lambda$-statistically quasi-Cauchy sequences of double sequences, and $\lambda$-statistically ward double continuity to find out whether $\lambda$-statistically ward double continuity coincides with $\lambda$- statistically ward (single) continuity or not (see \cite{PattersonandSavasLacunarystatisticalconvergenceofdoublesequences}, and \cite{SavasDoublelacunarystatisticallyconvergentsequencesintopologicalgroups} for the definitions and related concepts in the double case).

\end{document}